\documentclass{article}[12pt]
\usepackage{amsfonts, amsmath, theorem }
\newtheorem{lemma}{Lemma}[section]

\newtheorem{theorem}[lemma]{Theorem}
\newtheorem{corollary}[lemma]{Corollary}
{\theorembodyfont{\upshape}}
{\theorembodyfont{\upshape}
{\theorembodyfont{\upshape}}
{\theorembodyfont{\upshape}}
{\theorembodyfont{\upshape}}

\newcommand{\ep}{\varepsilon}
\newcommand{\Proof}{\underbar{Proof}{\hskip 0.1in}}
\newcommand{\rn}{\mathbb{R}^N}

\newcommand{\qt}{\frac{1}{4}}

\newcommand{\bd}{\partial \Omega}

\begin{document}
\title{PERTURBATION OF DOMAIN: ORDINARY DIFFERENTIAL EQUATIONS} 
\author{C. Mason}
\date{May 2001}
\maketitle

\begin{abstract}
We study a boundary perturbation problem for a one dimensional Schr\"odinger equation in which the potential has a regular singularity near the perturbed end point.  We give the asymptotic behaviour of the eigenvalues under the perturbation.  This problem arose out of the author's studies of singular elliptic operators in higher dimensions and we illustrate this point with an example.  The class of potentials to which this method applies is larger than that covered by standard results which assume uniform ellipticity of the operator or a perturbative term which is analytic in the perturbation parameter.  
\end{abstract}

\section{Introduction}

This paper forms part of the author's studies of the behaviour of the spectrum of the Dirichlet Laplacian under perturbations of the domain.  We let $\Omega$ be some bounded domain (either in $\rn$ or more generally in an incomplete Riemannian manifold) and then consider a slightly smaller region $\Omega_\ep$ such that
\[
\{ x\in\Omega : \text{dist}(x,\bd)<\ep \}\subseteq \Omega_\ep \subseteq \Omega
\]
for $\ep>0$ small enough.  Let $H$ denote an elliptic operator acting in $L^2(\Omega)$ with Dirichlet boundary conditions and suppose it has eigenvalues $\lambda_n$.  Now let $\lambda_{n,\ep}$ denote the eigenvalues of $H$ restricted to $\Omega_\ep$.  The problem then is to study the asymptotics of the quantity $ \lambda_{n,\ep}-\lambda_n$ as $\ep \to 0^+$. 

The specific model we shall consider in this paper is a one-dimensional Schr\"odinger equation
\begin{gather*}
-f''(x)+V(x)f(x)=\lambda f(x) \\
a<x<b
\end{gather*}
subject to Dirichlet boundary conditions.  The potential is assumed to have a regular singularity at the end point $x=a$ and thus we are led to consider the boundary perturbation problem with the perturbed domain $(a+\ep,b)$. More precisely we assume that
\[
V(x)(x-a)^2/c \to 1
\]
as $x\to a^+$ for some $-1/4<c<3/4$.  We do not however assume that $V(x)(x-a)^2$ is analytic and so we cannot employ standard techniques to determine series expansions of the eigenfunctions as described in \cite{O} for example.   

The stability of the spectrum under perturbations of the boundary has a long history.  In particular we mention the texts of Kato, \cite[Section 6.5]{K}, Courant-Hilbert, \cite[p419]{CH}, and Reed and Simon, \cite[Chapter XII]{RS4}.  These may be used to show that $\lambda_{n,\ep}$ are holomorphic in $\ep$ near $\ep=0$ provided the operator is uniformly elliptic and $\Omega_\ep$ may be obtained from $\Omega$ by a sufficiently well-behaved transformation.  These techniques work by changing the problem from one in which the domain is a function of some small parameter $\ep$ and the operator is fixed to one in which the domain is fixed and the operator depends on the small parameter in some analytic way.  

If the Euclidean region is star-shaped  then the eigenvalues of the Laplacian are known to satisfy
\[
\lambda_{n,\ep}-\lambda_n=O(\ep) \text{ as } \ep \to 0^+.
\]
However, in \cite{M1} the author considered a completely general type of boundary perturbation assuming no regularity of the coefficients or indeed of the boundary.  It merely assumed that some generalised Hardy inequality was valid.  We refer to this paper for details but stress that the results can be applied when the operator is singular, including operators of the form
\begin{align}\label{eqn:ode}
Hf:=-d(x)^{\alpha}\Delta_E
\end{align}
where $d(x)$ is the distance to the boundary, $0\leq \alpha<1$ and $\Delta_E$ is the ordinary Euclidean Laplacian in $\mathbb{R}^2$.  To show the optimality of the results we considered the operator (\ref{eqn:ode}) acting in a rotationally invariant domain and then used separation of variables to reduce the problem to the type of one-dimensional problem considered in this paper.  This specific problem is dealt with in Example \ref{ex:motiv}.  

The main result, contained in Theorem \ref{thm:main}, is that we have the asymptotic expansion
\[
\lambda_{n,\ep}=\lambda_n+c_n\ep^p+o(\ep^p)
\]
as $\ep \to 0^+$ where $p:=2\sqrt{c+1/4}$.  This result is important because it shows how to perform a boundary perturbation on a non-uniformly elliptic equation using only the standard theory of ODEs.  In particular the standard results of Kato cannot possibly be applied.  Moreover, Example \ref{ex:motiv} and its connection with (\ref{eqn:ode}) show that the result has implications for higher dimensional singular problems.

There are many good introductions to the theory of ordinary differential equations and in particular the classification of the boundary.  We will not reproduce this theory here.  The reader who is unfamiliar with the terminology or needs a refresher course could look for example at the comprehensive review of the subject in \cite{Z}.

Finally to fix notation we remark that we write
\[
f(x)\sim  g(x) \text{ as } x\to a+
\]
to mean
\[
\frac{f(x)}{g(x)} \to 1 \text{ as } x \to a+.
\]

\section{Set-up of the Problem}

Let $L$ be  the operator acting in $L^2(a,b)$  given initially by
\[
Lf:=-f''(x)+V(x)f(x) \text{ for } f\in C^2_c(a,b).
\]
We assume the following:
\begin{enumerate}
\item $V$ is continuous throughout $(a,b)$.
\item $V$ is twice differentiable near $x=a$.
\item $V$ satisfies
\[
V(x)\sim \frac{c}{(x-a)^2}\text{ as } x\rightarrow a+.
\]
where $-1/4<c<3/4$.  Moreover, this behaviour may be differentiated in the sense that $|V'|=O((x-a)^{-3})$ and $|V''|=O((x-a)^{-4})$. 
\item $L$ is bounded below and its spectrum consists solely of simple eigenvalues.
\end{enumerate}

Condition (1) will ensure that initial value problem is well-posed and has solutions in the classical sense.  Conditions (2) and (3) specify the singularity at the end point is regular and will allow us to use Liouville-Green theory to determine the asymptotic behaviour of the eigenfunctions.  The range of values of $c$ is not  arbitrary.  If $c<1/4$ then the corresponding operator ceases to be bounded below and if $c>3/4$ then the spectral properties change (see below).    

The simplicity of the spectrum is guaranteed by taking separated boundary conditions.  The fact that $L$ is bounded below allows us to take the Friedrichs extension and so deal with a self-adjoint operator.    We comment that the Hardy Inequality immediately gives us the fact that $L$ is bounded below provided $V_{-}$, the negative part of the potential, satisfies
\[
V_{-}(x)\geq -\frac{1}{4d(x)}
\]
for all $x\in(a,b)$ where $d(x)=\min \{|x-a|,|x-b|\}$ (condition (3) does not, of course, imply this inequality). 

\begin{lemma}\label{lem:approx}
Suppose 
\[
V(x) \sim \frac{c}{(x-a)^2} \text{ as } x \to a+
\]
where $c>-1/4$.  Moreover, suppose that this asymptotic behaviour can be differentiated. Then we have a basis of solutions comprising $f_+$ and $f_-$ where
\[
f_{\pm}(x)\sim(x-a)^{\frac{1}{2}\pm \sqrt{c+\frac{1}{4}}} \text{ as } x \to a+.
\]
Moreover, if $-1/4<c<3/4$ then both $f_+$ and $f_-$ belong to $L^2(a,\delta)$ for some $\delta>a$. 
\end{lemma}

\Proof  The behaviour of the eigenfunctions is proved using Liouville Green theory.  See \cite[Chapter 6, Section 3]{O} for details.  The convergence or otherwise of $\int_a^d |f_\pm(x)|^2dx$ is determined by the convergence or otherwise of $\int_a^d (x-a)^{1\pm 2\sqrt{c+\frac{1}{4}} }dx$ for $d$ near enough to $a$.

The solution $f_+$ is unique up to constant multiples and is refered to as the {\em recessive} solution.  The solution $f_2$ is non-unique and is called the {\em dominant} solution.

In the case that both $f_+$ and $f_-$ are square integrable near the end point and vanish only finitely often we must give a boundary condition (this is known as the `Limit Circle Non Oscillatory' case, denoted LCNO).  The problem of specifying the correct boundary condition so that the equation does in fact correspond to taking the Friedrichs extension of the operator $L$ has been solved.  We briefly outline the result here; for more information see \cite{Z}.

The maximal domain for $L$ on $(a,b)$, denoted by $M(a,b)$ is defined by
\begin{align*}
M(a,b):=&\left\{ f:(a,b)\to\mathbb{C} : f,f'\in AC(a,b) \text{ and } \right. \\ &\left. f,Lf \in L^2(a,b) \right\}
\end{align*} 
where $AC(a,b)$ denotes the space of absolutely continuous functions defined on $(a,b)$.
We define the Wronskian of  $u$,$v\in M(a,\delta]$ by
\[
[u,v](x):=-p(x)u'(x)v(x)+u(x)p(x)v'(x)
\]
and we make the definition
\[
[u,v](a):=\lim_{x \rightarrow a}[u,v](x).
\]
The correct LCNO Friedrichs boundary condition is then
\[
[f,u](a)=0
\]
where $u$ is a solution that satisfies $(u^2)^{-1}\not \in L^2(a,c)$ ($u$ is called a {\em principal solution at a} and in the non-oscillatory case is just the recessive solution).
\begin{corollary}\label{cor:bound}
Let $f_+$ be the recessive solution near $x=a$.  Then, in the case $-1/4<c<3/4$ the Friedrichs boundary condition is given by
\[
[f,f_+](a)=0
\]
\end{corollary}

If the constant $c$ is equal to $-1/4$ then Liouville Green cannot be used.  If it is the case that $(x-a)^2V(x)$ is analytic then one can use the Frobenius method to determine the solutions.  This procedure is outlined in \cite{O} but we do not consider this case in this paper.

\section{Boundary Perturbation}\label{sect:main}

\begin{theorem}\label{thm:main}
Let $L$ satisfy the assumptions above and denote the eigenvalues of $L$ by $\lambda_n$.  Now, let $L_\ep:=L\mid_{(a+\ep,b)}$ with Friedrichs boundary conditions and denote its eigenvalues by $\lambda_{n,\ep}$.  Then the following holds:
\[
\lambda_{n,\ep}=\lambda_n+c_n\ep^{p}+o(\ep^{p}) \text{ as } \ep \to 0+
\]
where
\[
p:=2\sqrt{c+1/4}.
\]
\end{theorem}

\Proof   Let $\lambda$ be an eigenvalue of $L$.  Near $x=a$ we have two linearly independent solutions, $\phi_{1,\lambda}$ and $\phi_{2,\lambda}$.  Since $x=a$ is LCNO both $\phi_{1,\lambda}$ and $\phi_{2,\lambda}$ lie in $L^2$ near $a$.  We suppose that $\phi_{1,\lambda}$ is the recessive solution and so satisfies the boundary condition at $x=a$ -- see Corollary \ref{cor:bound}.  Also,  we normalise the functions in order to assume that
\begin{align*}
\phi_{1,\lambda}(x)&\sim (x-a)^{1/2+\sqrt{c+1/4}} \\
\phi_{2  ,\lambda}(x)&\sim (x-a)^{1/2-\sqrt{c+1/4}}
\end{align*}
as $x\to a^+$.

Since $\lambda$ is an eigenvalue $\phi_{1,\lambda}$ also satisfies the boundary condition at $x=b$ (either an explicit condition if $b$ is R, WR or LCNO or the impicit condition that it must lie in $M(c,b)$ if $b$ is LP).

Now, we consider the differential equation
\[
-\xi_\lambda''(x)+(V(x)-\lambda)\xi_\lambda(x)=\phi_{1,\lambda}(x).
\]
At $x=a$ we impose the boundary condition
\[
[\xi_\lambda, \phi_{2,\lambda}](a)=0
\]
and at $x=b$ we use the same boundary condition as for the homogeneous equation.

The corresponding homogeneous equation is solved for the left boundary condition by $\phi_{2,\lambda}$ and for the right by $\phi_{1,\lambda}$.  These are linearly independent and so we can write down the Green's function:
\[
G(x,y):=\left\{ \begin{array}{cc} -\frac{1}{\kappa}\phi_{2,\lambda}(y)\phi_{1,\lambda}(x) & x \geq y \\
  -\frac{1}{\kappa}\phi_{1,\lambda}(y)\phi_{2,\lambda}(x) & y \geq x \end{array} \right.
\]
where $\kappa:=[\phi_{1,\lambda},\phi_{2,\lambda}](c)$ for any $c\in(a,b)$.  Note that $\kappa \ne 0$ and is a constant whose value is independent of the particular $c$ chosen -- see \cite[p351]{CH} for an introduction to Green's functions.

We write down the function
\[
h(x)=\phi_{1,\lambda}(x)+\ep^pc_\ep\xi_\lambda(x)
\]
where $c_\ep$ is some constant to be determined.  We know that this satisfies the right boundary condition and next turn to the point $x=a+\ep$.  We have
\[
h(a+\ep)=\phi_{1,\lambda}(a+\ep)+\ep^pc_\ep\xi_\lambda(a+\ep)
\]
and thus choosing
\begin{align*}
c_\ep&:=\frac{-\ep^{-p}\phi_{1,\lambda}(a+\ep)}{\xi_\lambda(a+\ep)} \\
&=\frac{-\ep^{-p}\phi_{1,\lambda}(a+\ep)}{-\frac{1}{\kappa}\int_a^{a+\ep}\phi_{2,\lambda}(y)\phi_{1,\lambda}(a+\ep )\phi_{1,\lambda}(y)dy-\frac{1}{\kappa}\int_{a+\ep}^{b} \phi_{1,\lambda}(y)\phi_{2,\lambda}(a+\ep)\phi_{1,\lambda}(y)dy}
\end{align*}
we also have
\[
h(a+\ep)=0.
\]
Using the asymptotic expansions developed for $\phi_{1,\lambda}$ and $\phi_{2,\lambda}$ we have the following:
\begin{align*}
\int_{a}^{a+\ep} \phi_{2,\lambda}(y)\phi_{1,\lambda}(y)dy=o(\ep^2) \text{ as } \ep \to 0+, \\
\int_{a+\ep}^{b}\phi_{1,\lambda}(y)\phi_{1,\lambda}(y)dy=\| \phi_{1,\lambda}\|_2^2+o(\ep^{2+p}) \text{ as } \ep \to 0^+.
\end{align*}
Thus we have
\begin{align*}
c_\ep&=\frac{\kappa \ep^{-p}\phi_{1,\lambda}(a+\ep)}{\phi_{2,\lambda}(a+\ep)(\| \phi_{1,\lambda}\|_2^2+o(\ep^{2+p}))} \\
&=\frac{\kappa}{\| \phi_{1,\lambda}\|_2^2}\frac{\ep^{-p}(\ep^{1/2+p/2}+o(\ep^{1/2+p/2})}{ (\ep^{1/2-p/2}+o(\ep^{1/2-p/2}))  }\frac{1}{(1+o(\ep^{2+p}))} \\
&=\frac{\kappa}{\| \phi_{1,\lambda}\|_2^2}+o(1) \text{ as } \ep \to 0^+.
\end{align*}

It is crucial to notice that the leading term in the expansion of $c_\ep$ is independent of the particular $\phi_{1,\lambda}$ and $\phi_{2,\lambda}$ chosen; i.e. suppose we now took $\widetilde{\phi_{1,\lambda}}=A\phi_{1,\lambda}$ and $\widetilde{\phi_{2,\lambda}}=B\phi_{2,\lambda}$, then a calculation shows that
\[
\widetilde{c_\ep}=c_\ep
\]
in an obvious notation.

Define $\lambda_\ep$ by
\[
\lambda_\ep:=\lambda+c_\ep \ep^p
\]
and the function $h$ satisfies
\[
-h''(x)+(V(x)-\lambda_\ep)h(x)=-c_\ep^2\ep^{2p}\xi_\lambda(x)
\]
subject to Friedrichs boundary conditions.

Now suppose $\lambda_\ep$ is not in the spectrum of $L_\ep$.  The operator $L_\ep-\lambda_\ep$ is thus invertible and we have
\[
\| (L_\ep-\lambda_\ep)^{-1}\xi_\lambda \|_2=\frac{\|h\|_2}{c_\ep^2\ep^{2p}}
\]
and so
\begin{align}\label{eqn:disttospec}
\| (L_\ep-\lambda_\ep)^{-1}\|\geq \frac{\| (L_\ep-\lambda_\ep)^{-1}\xi_\lambda \|_2}{\| \xi_\lambda \|_2}=\frac{\|h\|_2}{ \| \xi_\lambda \|_2  c_\ep^2\ep^{2p}}.
\end{align}
Let $C=\| h\|_2/\| \xi_\lambda \|_2$.  Then
\begin{align*}
C^{-1}=\frac{ \| \xi_\lambda \|_2}{\| h\|_2}&=\frac{\|\xi_\lambda \|_2}{\| \phi_{1,\lambda}+\ep^pc_\ep \xi_\lambda\|_2} \\
& \leq \frac{\|\xi_\lambda \|_2}{\|\phi_{1,\lambda}\|_2-c_\ep\ep^p\|\xi_\lambda\|_2} \leq 2\frac{\|\xi_\lambda \|_2}{\|\phi_{1,\lambda}\|_2}
\end{align*}
for $\ep$ small enough and thus we have a uniform upper bound for $C^{-1}$.

From (\ref{eqn:disttospec})  we immediately have
\begin{align*}
\frac{1}{\text{dist}(\lambda_\ep,\text{Spec}(L_\ep))}\geq \frac{C}{c_\ep^2\ep^{2p}} \\
\Leftrightarrow \text{dist}(\lambda_\ep,\text{Spec}(L_\ep)) \leq \frac{c_\ep^2\ep^{2p}}{C}.
\end{align*}

Thus there exists a point $\mu \in$ Spec$(L_\ep)$ such that
\begin{align*}
\mu &= \lambda_\ep+O(\ep^{2p}) \\
&=\lambda+c_\ep\ep^{p}+O(\ep^{2p})\\
&=\lambda+\left(\frac{\kappa}{\|\phi_{1,\lambda}\|_2}+o(1)\right)\ep^p+O(\ep^{2p}) \\
&=\lambda+\frac{\kappa}{\|\phi_{1,\lambda}\|_2}\ep^p+o(\ep^p) \text{ as } \ep \to 0+.
\end{align*}

This conclusion obviously still holds if $\lambda_\ep\in$ Spec($L_\ep$). 

We consider the bottom eigenvalue $\lambda_0$.  The above proves that there is a point $\mu \in$ Spec$(L_\ep)$ such that 
\[
\mu-\lambda_0 \leq c\ep^{p}
\]
for all $\ep$ small enough.  Thus $\mu$ eventually becomes smaller that $\lambda_1$ and hence smaller than $\lambda_{n,\ep}$ for all $n\geq 1$.  In other words $\mu=\lambda_{0,\ep}$ eventually.  The eigenvalues $\lambda_n$ all have multiplicity 1 and this argument can be repeated for $\lambda_1$ and so on.  This completes the proof.

\section{Examples}

\subsection{Example 1}
The theorem applied to the Bessel equation
\[
-f''(x)+\frac{\nu^2-1/4}{x^2}f(x)=\lambda f(x) \text{ for } 0<x<1
\]
in the case
\[
0<\nu<1
\]
yields that fact that
\begin{align}\label{eqn:asympt}
\lambda_{n,\ep}=\lambda_n+c_n\ep^{2\nu}+o(\ep^{2\nu}) \text{ as } \ep \to 0+.
\end{align}
This captures the main flavour of the results.  In particular it highlights the fact that the rate of convergence lies in the interval $(0,2)$.

\subsection{Example 2}\label{ex:motiv}

We now return to the main motivating example.  The original problem was two-dimensional.  We consider the operator
\[
Hf:=-d(x)^{2\gamma}\Delta_E
\]
acting the space $L^2(\Omega, d(x)^{-2\gamma})$ where $\Omega$ is the unit disc in two dimensions.  The rotational invariance of the domain allows us to separate variables and after some transformations reduce the problem to a Schrodinger equation.  Here the particular potential is \[
V(x)=\frac{2\gamma-\gamma^2}{4(1-\gamma)^2}\frac{1}{(\alpha-x)^2}+\left(\nu^2-\qt\right)\frac{(1-\frac{x}{\alpha})^{2\gamma \alpha}}{(1-(1-\frac{x}{\alpha})^{\alpha})^2}
\]
for $x\in(0,\alpha)$.  The constants $\gamma$ and $\alpha$ lie in the intervals $[0,1/2)$ and $[1,2)$ respectively and are related by $\alpha:=1/(1-\gamma)$.  In this case
\[
V(x)= \frac{2\gamma-\gamma^2}{4(1-\gamma)}\frac{1}{(\alpha-x)^2}+(1-\gamma)^{2\gamma\alpha}(\nu^2-1/4)(\alpha-x)^{2\gamma \alpha}+O((\alpha-x)^{(2\gamma+1)\alpha})
\]
as $x\to \alpha-$.  The results of the paper are then used to show that
\[
\lambda_{n,\ep}=\lambda_n+c_n\ep^{\frac{1}{1-\gamma}}+o(\ep^{\frac{1}{1-\gamma}}) \text{ as } \ep \to 0+.
\]
This is in agreement with the abstract result (see \cite[Theorem 7.1]{M1}).

\subsection{Numerical Example}
We carried out some numerical tests using the program SLEIGN2 - see \cite{BEZ} for an introduction.
For simplicity we return to the case of the Bessel equation on $(0,1)$ and choose $\nu=0.6$.  The eigenvalues are easily computed using  standard computing packages as well as SLEIGN.  Next for selected values of $\ep$ we can use SLEIGN to compute the bottom eigenvalue on $(\ep,1)$.  The results are recorded below.
\[
\begin{array}{|c|c|} \hline
\ep & \lambda_\ep \\ \hline
0.1 & 12.6988324 \\ \hline
0.01 & 10.8878760 \\ \hline
0.001 & 10.7821955 \\ \hline
0.0001 & 10.7755528 \\ \hline
0.00001 & 10.7751337 \\ \hline
0.000001 & 10.7751073 \\ \hline
\end{array}
\]
SLEIGN's quoted tolerances are all smaller than $10^{-10}$.  The computed value for $\lambda_0$ was
\[
\lambda_0=10.7751055.
\]

Assuming that a relation of the form
\[
\lambda_{0,\ep}-\lambda_\ep=c\ep^\delta
\]
making a plot of $\log (\ep)$ against $\log( \lambda_{0,\ep}-\lambda_\ep)$  we can estimate the value of $\delta$ from the slope of the graph.  Performing a least squares fit on the data in the table yields the value
\[
\delta=1.204292602.
\]
Although this is not a decisive piece of numerical analysis it does provide a reasonable confirmation of equation (\ref{eqn:asympt}) with $\nu=0.6$.

\textbf{Acknowledgements}  I would like to thank Brian Davies for suggesting this problem and for his guidance and advice during the work.  I also acknowledge the support of the Engineering and Physical Sciences Research Council through a research studentship.  I would also explicitly like to acknowledge the work of those involved in the SLEIGN project whose program has been useful to me.

Department of Mathematics \\
King's College London \\
Strand \\
London \\ 
WC2R 2LS \\
U.K.

email: cmason@mth.kcl.ac.uk

\end{document}